\documentstyle[12pt]{article}

\begin{document}
\newcommand{\s}{{\bf s}}
\centerline{\bf CLASSES OF BANACH SPACES STABLE AND UNSTABLE }
\centerline{\bf WITH RESPECT  TO THE OPENING}

\centerline{\bf M.I.Ostrovskii\footnote[1]{Supported by T\"UB\.ITAK}}
\centerline{Mathematics Department, Bo\v gazi\c{c}i University}
\centerline{80815 Bebek, Istanbul, TURKEY}
\centerline{and}
\centerline{Mathematical Division, Institute for Low Temperature Physics}
\centerline{47 Lenin avenue, 310164 Kharkov, UKRAINE}

e-mail: mostrovskii@ilt.kharkov.ua

{\bf Abstract.} The paper is a complement to the survey:
M.I.Ostrovskii "To\-po\-lo\-gies on the set of all subspaces of a Banach
space and related questions of Banach space geometry", Quaestiones
Math. (to appear). It contains proofs of some results on stability
of properties of Banach spaces with respect to the geometric opening
stated in the survey without proofs.

Some results of the present paper are of independent interest,
in particular the description of a predual property of the Banach--Saks
property.\footnote[2]{1991 Mathematics Subject Classification 46B20}
\bigskip



The present paper is  a  complement  to  the  survey [O1].
It contains proofs of results stated in [O1] without proofs.

Let us recall some of the notions from [O1]
and
fix some notation.
Let $X$ be a Banach space.
The closed unit ball and the unit sphere of
a Banach space $X$ are denoted by $B(X)$ and $S(X)$ respectively.
Recall that the {\it density character} of a  topological  space $X$  is
defined to be the least cardinality of dense  subset  of $X$  and  is
denoted by dens$X$.
We shall denote the set of all  closed
subspaces   of $X$   by $G(X)$.

Suppose $Y,Z\in G(X)$.    Let
$$\Theta _0(Y,Z)=\sup \{{\rm dist}(y,Z):\ y\in S(Y)\}$$
$$\Omega _0(Y,Z)=\sup \{{\rm dist}(y,S(Z)):\ y\in S(Y)\}.$$

It  is
clear that $0\le \Theta _0(Y,Z)\le \Omega _0(Y,Z)\le 2\Theta _0(Y,Z)$.

The {\it opening} between $Y$  and $Z$
is defined to be
$$
\Theta (Y,Z)=\max \{\Theta _0(Y,Z),\Theta _0(Z,Y)\}.
$$

The {\it spherical opening} between $Y$ and $Z$ is defined to be
$$
\Omega (Y,Z)=\max \{\Omega _0(Y,Z),\Omega _0(Z,Y)\}.
$$

We shall not distinguish terms ``class of  Banach
spaces'' and ``property of Banach spaces''.

Let $P$ be a class of Banach spaces. Class $P$ is called {\it stable}  if
there exists a number $\alpha >0$ such that for every  Banach  space $X$  and
every $Y,Z\in G(X)$ the following implication holds
$$(Y\in P)\&(\Theta (Y,Z)<\alpha)\Rightarrow (Z\in P).$$
Class $P$  is  called
{\it extendedly stable} if there exists a number $\alpha >0$ such that
for  every
Banach space $X$ and every $Y,Z\in G(X)$ the following implication holds
$$(Y\in P)\&(\Theta _0(Z,Y)<\alpha)\Rightarrow(Z\in P).$$

One of the methods of establishing of unstability presented  in
[O1] is based on the following construction.

Let $Y$  and $Z$  be  Banach  spaces  and  let $T:S(Y)\to S(Z)$   and
$D:S(Y^*)\to S(Z^*)$ be some surjective mappings (we suppose here that such
mappings exist). Let us introduce  on  the  algebraic  sum $Y\oplus Z$  the
following seminorm:
$$
p(y,z)=\sup \{|y^*(y)-(Dy^*)(z)|:\ y^*\in S(Y^*)\}.
$$
Seminorm $p$ generates  norm  on  the  quotient  of $Y\oplus Z$  by  the
zero-space of $p$. Denote the completion of this normed space by $X$.

By properties of $D$ the spaces $Y$ and $Z$ are  isometric to
their natural images in $X$. It is not hard to verify (see [O1], 6.14)
that the spherical opening between the images of $Y$ and $Z$ in $X$ is not
greater than the following value:
$$
\sup \{|y^*(y)-(Dy^*)(Ty)|:\ y\in S(Y),\ y^*\in S(Y^*)\}.
\eqno{(1)}$$
Let  us  introduce  the  quantity $k(Y,Z)$  as  the  infimum  of
quantities  (1) taken over  all  surjective  mappings $T:S(Y)\to S(Z)$   and
$D:S(Y^*)\to S(Z^*)$. If such mappings (at least one of them) do not exist
 we let $k(Y,Z)=2.$

Let us observe that the value (1) is not  greater  than  2  for
every $Y, Z, T, D$. So $k(Y,Z)\le 2$ for every pair $(Y,Z)$.

Following M.I.Kadets [K] we introduce for every pair $(Y,Z)$  of
Banach spaces the value
$$
d_{\Omega }(Y,Z)=\inf _{X,U,V}\Omega (UY,VZ)
$$
where the infimum is taken  over  all  Banach  spaces $X$  containing
isometric copies of $Y$ and $Z$ and over all isometric embeddings $U:Y\to X$
and $V:Z\to X$.

Arguments above show that
$d_{\Omega }(Y,Z)\le k(Y,Z)$.
In [O1] I asserted that these values are equivalent.  Now  we  shall
prove this assertion.

{\bf Proposition 1.} {\it For every Banach spaces $Y$ and $Z$  the  following
inequality
$$
20d_{\Omega }(Y,Z)\ge k(Y,Z).
$$
is valid.}

Proof. From the inequality $k(Y,Z)\le 2$ it follows that we  may  restrict
ourselves to the case $d_{\Omega }(Y,Z)<1/10.$

Let $\varepsilon $ be an arbitrary number of the  open
interval  (0,1).  By
definition of $d_{\Omega }$ there  exists  a  Banach
space $X$
containing subspaces isometric to $Y$ and $Z$ (we  shall  still  denote
them by $Y$ and $Z)$, for which
$$
\Omega (Y,Z)\le d_{\Omega }(Y,Z)(1+\varepsilon ).
$$
We shall use the following definition. Subset $L$  of  a  metric
space $X$ is said to be $\tau$-{\it lattice} $(\tau\in {\bf R},\
\tau>0)$, if the distance  between each two elements of $L$ is not less
than $\tau$, and is said to be {\it maximal}
$\tau$-{\it lattice}, if $L$ is not a proper subset of
any $\tau$-lattice in $M$. By standard application of Hausdorff's
maximality theorem it follows that every metric space contains
a maximal $\tau$-lattice for arbitrary $\tau>0$.

Let $\{y_{\alpha }\}_{\alpha \in A}$ be a maximal
$(2(1+\varepsilon )\Omega (Y,Z))$-lattice  in $S(Y)$.  By
definition of maximal $\tau$-lattice we can  find  such
subsets $B_{\alpha }\subset S(Y)\ (\alpha \in A)$ that
the following conditions are satisfied.

(a) $B_{\alpha }$ contains the intersection of $S(Y)$ with the open
ball of radius $(1+\varepsilon )\Omega (Y,Z)$ centered  at $y_{\alpha }$.

(b) $B_{\alpha }$ is contained in the ball  of the radius
$2(1+\varepsilon )\Omega (Y,Z)$
centered at $y_{\alpha }$.

(c) Sets $B_{\alpha }$ are pairwise disjoint.

(d) $\cup _{\alpha \in A}B_{\alpha }=S(Y).$

For every $\alpha \in A$ we choose $z_{\alpha }\in S(Z)$ such that
$||z_{\alpha }-y_{\alpha }||<(1+\varepsilon )\Omega (Y,Z)$.
It is easy to verify that $\{z_{\alpha }\}_{\alpha \in A}$ is a
$4(1+\varepsilon )\Omega (Y,Z)$-net in $S(Z)$. Let
us  denote  the  intersection  of $S(Z)$  with  the  ball  of  radius
$4(1+\varepsilon )\Omega (Y,Z)$ centered at $z_{\alpha }$ by
$R_{\alpha }\ (\alpha \in A)$. We have $S(Z)=\cup _{\alpha \in A}R_{\alpha }$.

Since we may (and shall) assume without loss of generality that
$\Omega (Y,Z)<1/2$ (see remark at the begining of the  proof),  then  dens$Y=
$dens$Z$ (see [KKM] (Theorem 3) or [O1] (section  6.27.2)).  Therefore
by condition (a) card$R_{\alpha }={\rm card}B_{\alpha }$.
So we can  define  a
mapping $T:S(Y)\to S(Z)$ in such a way that $Ty_{\alpha }=z_{\alpha }$
and $T(B_{\alpha })=R_{\alpha }$.  It  is
clear that this mapping is surjective.

Let $y^*\in S(Y^*)$. It is easy to see that,  if
we extend it in a norm-preserving manner onto the whole $X$  and  then
restrict this extension to $Z$, then the norm of the obtained
functional is  not
less than $1-\Omega _0(Y,Z)$.

Let $\{y^*_{\gamma }\}_{\gamma \in C}$ be some maximal
$(2(1+\varepsilon )\Omega (Y,Z))$-lattice in $S(Y^*)$. It
is easy to see that there exists subsets
$B^*_{\gamma }\subset S(Y^*)\ (\gamma \in C)$  satisfying
the following conditions:

  (a) $B^*_{\gamma }$ contains the intersection of $S(Y^*)$  with  the  open  ba
ll  of
radius $(1+\varepsilon )\Omega (Y,Z)$ centered at $y^*_{\gamma }$.

  (b) $B^*_{\gamma }$ is contained in the ball  of  radius
$2(1+\varepsilon )\Omega (Y,Z)$ centered at $y^*_{\gamma }$.

  (c) Sets $B^*_{\gamma }$ are pairwise disjoint.

(d) $\cup _{\gamma \in C}B^*_{\gamma }=S(Y^*).$

Let $\{f_{\gamma }\}_{\gamma \in C}$ be some
norm-preserving extensions  of  functionals
$\{y^*_{\gamma }\}_{\gamma \in C}$ onto $X$.
Let $z^*_{\gamma }=f_{\gamma }|_Z/||f_{\gamma }|_Z||$.
Show that $\{z^*_{\gamma }\}_{\gamma \in C}$  is  an
$(10(1+\varepsilon )\Omega (Y,Z))$-net in $S(Z^*)$.

Suppose the contrary. Then  there  exists $z^*\in S(Z^*)$  such  that
$||z^*-z^*_{\gamma }||>10(1+\varepsilon )\Omega (Y,Z)$ for every $\gamma \in $C.

Therefore for every $\gamma \in C$  there
exists $z_{\gamma }\in S(Z)$ such that
$$
|(z^*-f_{\gamma }/||f_{\gamma }|_Z||)(z_{\gamma })|>10(1+\varepsilon )\Omega (Y,
Z).
$$

Let $g$ be some norm-preserving  extension  of $z^*$  onto  $X$.  Let
$y_{\gamma }\in S(Y)\ (\gamma \in C)$ be such that
$||z_{\gamma }-y_{\gamma }||<(1+\varepsilon )\Omega (Y,Z)$. We have
$$
|(g-f_{\gamma }/||f_{\gamma }|_Z||)(y_{\gamma })|>8(1+\varepsilon )\Omega (Y,Z).
$$
Therefore
$$
||||f_{\gamma }|_Z||g|_Y-y^*_{\gamma }||>(1-\Omega (Y,Z))8(1+\varepsilon )\Omega
 (Y,Z).
$$
Since we assumed that $\Omega (Y,Z)<1/2,$ then
$$
||||f_{\gamma }|_Z||g|_Y-y^*_{\gamma }||>4(1+\varepsilon )\Omega (Y,Z).
$$
Since the lattice $\{y^*_{\gamma }\}_{\gamma \in C}$ is maximal,
it follows that for  some
$\gamma \in C$ we have
$$
||g|_Y/||g|_Y||-y^*_{\gamma }||<2(1+\varepsilon )\Omega (Y,Z).
$$
Hence
$$
||y^*_{\gamma }-||f_{\gamma }|_Z||g|_Y||<
2(1+\varepsilon )\Omega (Y,Z)+||g|_Y||
(1/||g|_Y||-||f_{\gamma }|_Z||)<4(1+\varepsilon )\Omega (Y,Z).
$$
This contradiction proves the assertion concerning
$\{z^*_{\gamma }\}_{\gamma \in C}$.  Let  us
denote  the  intersection  of $S(Z^*)$  with  the   ball   of   radius
$10(1+\varepsilon )\Omega (Y,Z)$ with the centre in
$z^*_{\gamma }$ by $R^*_{\gamma }\ (\gamma \in C)$. It  is  clear  that
$S(Z^*)=\cup _{\gamma \in C}R^*_{\gamma }$.

By $\Omega (Y,Z)<1/2$ it follows (see e.g. [O2])  that  dens$Y^*$=dens$Z^*$.
Therefore for every $\gamma \in C$ we have card$R^*_{\gamma }$=card$
B^*_{\gamma }$. Hence we can  find  a
mapping $D:S(Y^*)\to S(Z^*)$ such that
$D(y^*_{\gamma })=z^*_{\gamma }$ and
$D(B^*_{\gamma })=R^*_{\gamma }$ for every $\gamma \in C$.
It is clear that $D$ is surjective.

Let us prove that the mappings $D$ and $T$  give us  the
desired estimate. Let $y^*\in B^*_{\gamma }$ and $y\in B_{\alpha }$. We have
$$|y^*(y)-(Dy^*)(Ty)|\le$$
$$
|y^*_{\gamma }(y_{\alpha })-y^*_{\gamma }(y_{\alpha }-y)-(y^*_{\gamma }-y^*)y-(z
^*_{\gamma }(z_{\alpha })-z^*_{\gamma }(z_{\alpha }-Ty)-(z^*_{\gamma }-Dy^*)(Ty)
)|\le
$$
$$
|y^*_{\gamma }(y_{\alpha })-z^*_{\gamma }(z_{\alpha })|+2(1+\varepsilon )\Omega
(Y,Z)+2(1+\varepsilon )\Omega (Y,Z)+
$$
$$
4(1+\varepsilon )\Omega (Y,Z)+10(1+\varepsilon )\Omega (Y,Z)\le
$$
$$
|f_{\gamma }(y_{\alpha }-z_{\alpha })-(f_{\gamma }-z^*_{\gamma })(z_{\alpha })|+
18(1+\varepsilon )\Omega (Y,Z)\le (1+\varepsilon )\Omega (Y,Z)+\Omega (Y,Z)+
$$
$$
18(1+\varepsilon )\Omega (Y,Z)\le 20(1+\varepsilon )\Omega (Y,Z).
$$
Hence
$$
k(Y,Z)\le 20(1+\varepsilon )\Omega (Y,Z)\le
20d_{\Omega }(Y,Z)(1+\varepsilon )^{2}.
$$
Since $\varepsilon $ is arbitrary number from the open  interval  (0,1),  we
obtain the desired inequality. The proposition is proved.$\Box$

The aim of the remaining part of the present paper  is  to  support
the conjecture on  characterization  of  extendedly  stable  classes
which was formulated in [O1] (section 6.36).

Recall necessary definitions.

Let $\Gamma $ be a set, and let $l_1(\Gamma )$ be  the  corresponding  space  of
functions on $\Gamma $. Let $X$ be a Banach space. By $l_{\infty }(\Gamma ,X)$
we  denote  the
space of bounded functions $x:\Gamma \to X$ with the $\sup $-norm.

Let $A$ be a subset of the unit sphere of $l_1(\Gamma )$.
For every $a\in A$  we
introduce a linear operator from $l_{\infty }(\Gamma ,X)$  into $X$  defined  by
  the
rule
$$
x\to \sum_{\gamma \in \Gamma }a(\gamma )x(\gamma ).
$$
This operator will be also denoted by $a$. It is clear that  the  norm
of this operator equals 1. The supremum of those
$\delta $ for  which  there
exists $x\in S(l_{\infty }(\Gamma ,X))$ such that
$\inf _{a\in A}||a(x)||\ge \delta $ is called the {\it index  of
$A$ in} $X$ and is denoted by $h(X,A)$

It is known (see [O1] (section 6.27)) that for many common properties
of Banach spaces (reflexivity, $B$-convexity, Banach-Saks property)
sets $\Gamma,\ A\subset S(l_1(\Gamma ))$
and $\tau>0$ can be chosen in such a  way  that
the absence of the property for a Banach space $X$  is  equivalent  to
each of the following inequalities:
$$
h(X,A)>0;
$$
$$
h(X,A)\ge \tau.
$$
Properties   admitting    such    descriptions    are    called
$l_1${\it -properties.} It is not  hard  to  verify
that $l_1$-properties  are
extendedly stable (see [O1], Proposition 6.21).

In [O1] I conjectured that the  converse  is  also  true,  more
precisely, that for every cardinal number $\alpha $ the intersection of
every extendedly stable property with the set  of  Banach  spaces  with
density character less than $\alpha $ is an $l_1$-property.

We shall prove several results supporting this  conjecture.  In
passing we shall find new classes of $l_1$-properties.

We introduce a subclass of the class of all $l_1$-properties  and
prove that all known at the moment extendedly stable classes  belong
to this subclass.

{\bf Definition 1.} Class $P$ of Banach spaces will be called a {\it regular}
$l_1${\it -property} if there  exist a real  number $\delta >0,$
a set $\Gamma $  and a subset
$A\subset S(l_1(\Gamma ))$ satisfying the conditions:

1. The set $A$ consists of finitely non-zero vectors.

2. If $a_0\in A$ then $A$ contains all vectors
$a\in S(l_1(\Gamma ))$ for which
$$
(\forall \gamma \in \Gamma )(\hbox{sign}a_0(\gamma )=\hbox{sign}a(\gamma )).
$$

3. For a Banach space $X$ the following conditions are equivalent:

(a) $X\not\in P.$

(b) $h(X,A)>0.$

(c) $h(X,A)\ge \delta .$

Supremum of those $\delta >0$ for which there exist $\Gamma $
and $A\subset S(l_1(\Gamma ))$  such
that conditions 1, 2 and 3 are satisfied will be called the  {\it regular
exponent} of class $P$ and will be denoted by reg($P$).

When working with regular $l_1$-properties, we shall use the  following
definition and notation. Let $A\subset S(l_1(\Gamma ))$.
Mapping {\bf s}$:\Gamma \to \{0,1,-1\}$  will
be called an {\it admissible sign} for $A$ if there
exists $a\in A$ such that
$$
(\forall \gamma \in \Gamma )({\bf s}(\gamma )=\hbox{sign}a(\gamma ))
\eqno{(2)}$$

The set of all admissible signs for $A$ will be denoted by {\bf S}$(A)$.
Let {\bf s}  be  an  admissible
sign for $A$. The set of all $a\in A$ for which (2) is  satisfied  will  be
denoted by $A(${\bf s}). It is clear that
$A=\cup _{{\bf s}\in{\bf S}(A)}A({\bf s})$.

When working with regular $l_1$-properties we shall repeatedly use  the
following lemma.

{\bf Lemma 1.} {\it Let $\Gamma $ be arbitrary set and let subset
$A\subset S(l_1(\Gamma ))$ be such
that conditions {\rm 1} and {\rm 2} of Definition {\rm 1} are  satisfied.
Suppose the Banach
space $X$ and $x\in S(l_{\infty }(\Gamma ,X))$ are such that
$$
\inf _{a\in A}||a(x)||\ge c.
\eqno{(3)}$$
Then for every {\bf s}$\in${\bf S}$(A)$ there exists a functional
$f_{{\bf s}}\in B(X^*)$  such  that
for every $\gamma \in ${\rm supp}{\bf s} the following relations are valid:}
$$
|f_{{\bf s}}(x(\gamma ))\mid \ge c.$$
$$
\hbox{sign}f_{{\bf s}}(x(\gamma ))={\bf s}(\gamma ).
$$

Proof. Consider the set
$\{{\bf s}(\gamma )x(\gamma ):\ \gamma \in $supp{\bf s}$\}$. By condition
2 of Definition 1 and inequality (3), the convex hull
of this set does not intersect the open ball of radius $c$
centered at zero. Using  the  separation  theorem,  we  find  required
functional.$\Box$

Now we are going to show that all $l_1$-properties listed in [O1]
(section 6.27) are regular $l_1$-properties. Here I should note that at
the moment I do not know $l_1$-properties which are not regular.

{\bf 1.} Reflexivity. Description presented  in [O1] (in  6.19)  doesn't
satisfy condition 2 of  Definition  1.
Nevertheless reflexivity is a regular $l_1$-property.
In  order  to prove this assertion
we need another description of reflexivity  (see  [B1,
p.~51]):

Let $X$  be  a  Banach  space.  The  following  conditions   are
equivalent:

  (a) $X$ is non-reflexive.

  (b) For every $\theta ,\ 0<\theta <1,$
there exist a sequence $\{x_k\}_{k\in {\bf N}}$ in $S(X)$ and a
sequence $\{f_{n}\}_{n\in {\bf N}}$ in $S(X^*)$ such that:
\[ f_{n}(x_k)=\left\{\begin{array}{ll} \theta , & {\rm if\ } n\le k,\\
0, & {\rm if\ } n>k.\end{array}
\right. \]

Let $A\subset S(l_1)$  be  the  set  of  all  finitely  non-zero  vectors
$a=\{a_k\}^{\infty }_{k=1}$ such that for some $n\in {\bf N}$ (depending on $a$)
we have $a_k\le 0$  when
$k\ge n$ and $a_k\ge 0$  when $k<n$.
It  is  clear  that  this  set  satisfies
conditions 1 and 2 of Definition 1.

Let $X$ be a Banach  space.  Let  us  prove  that  the  following
conditions are equivalent.

(a) $X$ is non-reflexive.

(b) $h(X,A)>0.$

(c) $h(X,A)\ge 1/3.$

Since the set $A$ contains the set considered in [O1] (6.19),  we
need to prove implication (a)$\Rightarrow$(c) only.

Let $a\in A$ and let $n\in {\bf N}$ be such that
$a_k\le 0$ when $k\ge n$  and $a_k\ge 0$  when
$k<n$. We have
$$
f_1(\sum^{\infty }_{k=1}a_kx_k)=\theta (\sum^{n-1}_{k=1}|a_k|-
\sum^{\infty }_{k=n}|a_k|).
$$
$$
f_{n}(\sum^{\infty }_{k=1}a_kx_k)=-\theta \sum^{\infty }_{k=n}|a_k|.
$$
Hence
$$
(f_1-2f_{n})(\sum^{\infty }_{k=1}a_kx_k)=
\theta \sum^{\infty }_{k=1}|a_k|=\theta .
$$
Therefore $||a(x)||\ge \theta /3.$ Since the description
above is valid for every
$\theta <1,$ it follows that $h(X,A)\ge 1/3.$

{\bf 2.}
It is known (see [O1], Example 6.20) that  for  every
uncountable cardinal $\alpha $ the class of all Banach spaces  with  density
character less than $\alpha $ is an $l_1$-property. But the sets
introduced  in
[O1] do not satisfy Condition 2 of Definition 1. Nevertheless  these
properties are regular $l_1$-properties.  This  can  be  shown  in  the
following way. Let $\alpha $ be an uncountable cardinal and
let $\Gamma $ be  a  set
of cardinality $\alpha $. Let us introduce subset
$A\subset S(l_1(\Gamma ))$ as the  set  of
all vectors with two-point  support.  It  is  clear  that  this  set
satisfies conditions 1 and 2 of Definition 1.

Let $X$  be  a  Banach  space.  The  following  conditions   are
equivalent:

(a) dens$X\ge \alpha $.

(b) $h(X,A)>0.$

(c) $h(X,A)\ge 1/3.$

Since the introduced set contains one considered  in [O1] (6.20),
then we need to prove the implication (a)$\Rightarrow$(c) only.

Let $X$ be a Banach space with dens$X\ge \alpha $.
We may consider $\Gamma $  as  a
set of all ordinals  which  are  less  than  the  least  ordinal  of
cardinality $\alpha $. By well-known arguments (see [KKM,  p.  98])  we  can
find subset $\{x(\gamma )\}_{\gamma \in \Gamma }\subset S(X)$
such that for every $\gamma \in \Gamma $ we have
$$
\hbox{dist}(x(\gamma ),\hbox{lin}\{x(\beta ):\ \beta <\gamma \})>
1-\varepsilon .
$$
It is easy to verify that for every pair $\{\beta ,\gamma \}\subset \Gamma $
with $\beta \neq \gamma $  we
have
$$
\inf \{||rx(\gamma )+sx(\beta )||:\ r,s\in{\bf R},\ |r|+|s|=1\}\ge
(1-\varepsilon)/(3-\varepsilon).
$$
Therefore for every $a\in A$ we have $||a(x)||\ge(1-\varepsilon)/
(3-\varepsilon).$ Since $\varepsilon$ is arbitrary, then $h(X,A)\ge1/3$.

{\bf3.} The class of finite dimensional Banach spaces. In this case we
take $\Gamma={\bf N}$ and $A\subset S(l_1)$ being
the set of all vectors with
two-point support. Using arguments from the previous example it is easy
to verify that for a Banach space $X$ the following statements
are equivalent

(a) $X$ is infinite dimensional

(b) $h(X,A)>0.$

(c) $h(X,A)\ge 1/3.$

Furthermore, using the main result of [P] it can be shown that the statements
(a)--(c) are equivalent to $h(X,A)\ge 1/2.$

The sets introduced in [O1] (see 6.27) in order to describe  classes 4--7 below
satisfy all conditions of Definition 1. So these
classes are regular $l_1$-properties.

{\bf 4}. The class of Banach spaces whose dimension is not greater than
$n$ (for every $n\in{\bf N}$).

{\bf 5}. $B$-convexity.

{\bf 6}. The class of Banach spaces not containing isomorphic copies
of $l_1$.

{\bf 7}. Alternate signs Banach-Saks property.

{\bf 8}. Super-reflexivity. Subset of $S(l_1)$ introduced in
[O1] to describe this property does not satisfy condition 2 of
Definition 1. Therefore we introduce another subset $A\subset S(l_1)$.
Let $A$ be the set of all vectors of the following form:
$$(0,\dots,0,a_1,\dots,a_n,0,\dots),$$
where $\{a_i\}^n_{i=1}$ are such that for some
$k\in\{1,\dots,n+1\}$ we have $a_i\ge 0$ when $i<k$ and $a_i\le
0$ when $i\ge k$, and $a_1$ is preceded by $n(n-1)/2$ zeros.

It is clear that $A$ satisfies conditions 1 and 2 of Definition 1.

Let $X$  be  a  Banach  space. Prove that the  following
conditions   are
equivalent:

(a) $X$ is not super-reflexive.

(b) $h(X,A)>0.$

(c) $h(X,A)=1.$

Since the set $A$ contains the set considered in [O1] the only thing
which we need to prove is that (a) implies (c). In order to do this
let us recall [B1, p. 270--271] that in every non-super-ref\-le\-xive
space for every $\varepsilon>0$ and $n\in{\bf N}$ there exists a set
$\{x_1,x_2,\dots,x_n\}\subset B(X)$ such that for every $1\le k\le n$
we have
$$||x_1+\dots+x_k-x_{k+1}-\dots-x_n||\ge n(1-\varepsilon).$$

Let $\{a_i\}_{i=1}^n$ be any sequence satisfying the conditions above
with $k\le n$ (the case $k=n+1$ is easy). It is clear that
$$||\sum_{i=1}^na_ix_i||\ge ||x_1+\dots+x_k-x_{k+1}-\dots-x_n-
\sum_{i=1}^k(1-a_i)x_i+\sum_{i=k+1}^n(1+a_i)x_i||\ge$$
$$n(1-\varepsilon)-(\sum_{i=1}^k|1-a_i|+\sum_{i=k+1}^n|1+a_i|)=$$
$$n(1-\varepsilon)-\sum_{i=1}^k(1-|a_i|)-\sum_{i=k+1}^n(1-|a_i|)=$$
$$n(1-\varepsilon)-(n-1)=1-n\varepsilon.$$

The implication (a)$\Rightarrow$(c) can be easily derived from this
estimate.

{\bf 9}. Banach--Saks property.

Consider the set $A\subset S(l_1)$ consisting of all vectors
satisfying the conditions

(a) The cardinality of support is not greater than the least
of its elements.

(b) For every $\{a_k\}_{k=1}^\infty$ there exists $n\in{\bf N}$ such that
$a_k\ge 0$ if $k<n$ and $a_k\le0$ if $k\ge n$.

B.Beauzamy [B2] proved that a Banach space has Banach--Saks
property if and only if it is simultaneously reflexive and has
alternate signs Banach--Saks property. Using this result,
characterization of alternate signs
Banach--Saks property
(see [O1] (6.27.8))
and the above characterization of reflexivity
it is not hard to verify that
the following statements are equivalent:

(a) $X$ is not Banach--Saks

(b) $h(X,A)>0.$

(c) $h(X,A)\ge1/3.$

{\bf Proposition 2.} {\it Let $\{P_\alpha\}_{\alpha\in C}$ be some
set of regular
$l_1$-properties and the inequality
$$\inf_\alpha{\rm reg}P_\alpha>0$$
is satisfied. Then $\cup_{\alpha\in C}P_\alpha$ is a regular $l_1$-property
and} ${\rm reg}(\cup_{\alpha\in C}P_\alpha)\ge\inf_\alpha{\rm reg}P_\alpha$.

Proof. Let $\varepsilon$ be an arbitrary number from the open interval
(0,1). Let $\Gamma_\alpha$ and $A_\alpha\subset S(l_1(\Gamma_\alpha))\
(\alpha\in C)$ be such that $A_\alpha$ satisfies conditions 1 and 2 of
Definition 1 and for a Banach space $X$ the following statements are
equivalent:

(a) $X\notin P_\alpha.$

(b) $h(X,A_\alpha)>0.$

(c) $h(X,A_\alpha)>(1-\varepsilon)$reg$P_\alpha$.

Let $\Gamma=\cup_{\alpha\in C}\Gamma_\alpha$. Let $A\subset S(l_1(\Gamma))$
be the union of natural images of $A_\alpha$ in $S(l_1(\Gamma))$. It is clear
that $A$ satisfies conditions 1 and 2 of Definition 1.

Let $X$ be a Banach space. It is easy to see that the following statements
are equivalent:

(a) $X\notin (\cup_{\alpha\in C}P_\alpha)$.

(b) $h(X,A)>0$.

(c) $h(X,A)\ge(1-\varepsilon)\inf_{\alpha\in C}\hbox{reg}P_\alpha$.

Since $\varepsilon\in$(0,1) is arbitrary, the proposition is proved.
$\Box$

{\bf Definition 2.} Let $P$ be some class of Banach spaces. {\it Preclass}
of $P$ is defined to be the class of all Banach spaces whose duals
belong to $P$. Preclass of $P$ is denoted by pre$(P)$.

It is known that if $P$ is extendedly stable then pre($P)$ is also
extendedly stable (for different versions and proofs of this result
see [J], [AAG] and [O1] (Proposition 6.34)).

{\bf Theorem 1.} {\it If $P$ is a regular $l_1$-property then
{\rm pre}($P$) is also
a regular $l_1$-property and} reg(pre($P))\ge$reg($P$).

Proof. Let $\theta$ be an arbitrary number from the open interval
(0,1). Let us choose a set $\Gamma$ and a subset $A\subset S(l_1(\Gamma))$
in such a way that all conditions of Definition 1
with $\delta=\theta$reg($P$) are satisfied. Put $\Psi={\bf S}(A)$. Let
$D\subset S(l_1(\Psi))$ be a subset consisting of all finitely non-zero
vectors $d$ satisfying the following condition. If
supp$d=\{{\bf s}_1,\dots,{\bf s}_n\}$, then for some
$$\gamma\in\cap_{i=1}^n\hbox{supp \bf s}_n$$
we have
$$(\forall i\in\{1,\dots,n\})(\hbox{sign} d({\bf s}_i)={\bf s}_i(\gamma)).$$
It is clear that $D$ satisfies the conditions 1 and 2 of
Definition 1. In order to prove Theorem 1 it is sufficient to
prove that the following statements are equivalent.

(a) $X\notin$pre$(P)$.

(b) $h(X,D)>0$.

(c) $h(X,D)\ge\theta$reg$(P)$.

It is clear that we need to prove implications (a)$\Rightarrow$(c)
and (b)$\Rightarrow$(a) only.

Let $\delta$ be an arbitrary number from the open interval (0,1).
Let $X\notin$pre$(P)$. It means that $X^*\notin P.$ Let vector
$x^*\in S(l_\infty(\Gamma,X^*))$ be such that
$$\inf_{a\in A}||a(x^*)||>\theta\delta\hbox{reg}P.$$
By Lemma 1 we can find for every ${\bf s}\in{\bf S}(A)$ a functional
$f_{\bf s}\in B(X^{**})$ such that for every $\gamma\in\hbox{supp}{\bf s}$
we have $|f_{\bf s}(x^*(\gamma))|>\delta\theta{\rm reg}(P)$ and
${\rm sign}f_{\bf s}(x^*(\gamma))={\bf s}(\gamma)$. It is clear
that we may suppose that $||f_{\bf s}||<1$.

Further we shall repeatedly use the following statement which goes
back to E.Helly (see [B1, p. 52]). Let $f\in X^{**}$ and
$||f||<1$. Let $\{x^*_1,x^*_2,\dots,x^*_n\}$ be a finite subset
of $X^*$. Then there exists a vector $x\in X$ such that
$||x||<1$ and $x^*_i(x)=f(x^*_i)$ for every $i\in\{1,\dots,n\}$.
Every such vector we shall call a {\it reflection of $f$ with respect
to} $\{x^*_1,x^*_2,\dots,x^*_n\}$.

Let $\{x(\s)\}_{\s\in\Psi}$ be some reflections of $\{f_\s\}_{\s\in\Psi}$
with respect to sets
$$\{x^*(\gamma)\}_{\gamma\in{\rm supp}\s}.$$
Let
$d\in D$. Let us estimate from below the value $||d(x)||$. Let
supp$d=\{\s_1,\s_2,\dots,\s_n\}.$ Let $\gamma\in\cap_{i=1}^n$supp$\s_i$
be such that sign$d(\s_i)=\s_i(\gamma)$ for $i\in\{1,\dots,n\}.$ We have
$$||d(x)||=||\sum_{i=1}^nd(\s_i)x(\s_i)||\ge x^*(\gamma)
(\sum_{i=1}^nd(\s_i)x(\s_i))=$$
$$\sum_{i=1}^nd(\s_i)f_{\s_i}(x^*(\gamma))=
\sum_{i=1}^n\s_i(\gamma)|d(\s_i)|\s_i(\gamma)|f_{\s_i}(x^*(\gamma))|\ge
\delta\theta{\rm reg}(P).$$
Since $\delta$ is arbitrary number from the open interval (0,1),
the argument above implies that (a)$\Rightarrow$(c).

Let us turn to the implication (b)$\Rightarrow$(a). Let
$h(X,D)>\varepsilon>0$. By Lemma 1 and weak$^*$ compactness
of $B(X^*)$, for every $\gamma\in\Gamma$
we can find $x^*(\gamma)\in B(X^*)$ such that
$|x^*(\gamma)x(\s)|>\varepsilon$ and
sign$(x^*(\gamma)x(\s))=\s(\gamma)$ for every {\s} for which
$\gamma\in$supp\s. Let $a\in A(\s)$. We have
$$||\sum_{\gamma\in\Gamma} a(\gamma)x^*(\gamma)||\ge
(\sum_{\gamma\in{\rm supp}\s} a(\gamma)x^*(\gamma))x(\s)=
\sum_{\gamma\in\Gamma} |a(\gamma)|\s(\gamma)|x^*(\gamma)x(\s)|
\s(\gamma)>\varepsilon.$$
Hence $X^*\notin P$ and $X\notin{\rm pre}(P)$. Theorem 1 is
proved. $\Box$

{\bf Remark.} If property $P$ has a concrete description as a regular
$l_1$-property, then the proof of Theorem 1 gives us a concrete
description of pre$(P)$. Hence we obtain e.g. a description of
pre--Banach--Saks property.

{\bf Definition 3.} Let $P$ be a class of Banach spaces. {\it Coclass} of
$P$ is defined to be the class of all Banach spaces for which the
quotient space $X^{**}/X$ is in $P$. Coclass of $P$ is denoted by
$P^{\rm co}$.

J.Alvarez, T.Alvarez and M.Gonzalez [AAG] proved that the extended
stability of $P$ implies the extended stability of $P^{\rm co}$. Our
aim is to prove the following result.

{\bf Theorem 2.} {\it Let $P$ be a regular $l_1$-property. Then
$P^{\rm co}$ is also a regular $l_1$-property.}

Proof. Let a set $\Gamma$, a subset $A\subset S(l_1(\Gamma))$ and
a number $\delta>0$ be such that all conditions of Definition 1
are satisfied. Let $\Psi$ be the set of all triples
$(\gamma,\alpha,k)$, where $\gamma\in\Gamma$, $\alpha$ is a finite
subset of ${\bf S}(A)$, $k\in{\bf N}$.

We define subset $E\subset S(l_1(\Psi))$ in the following way.
Let $\{\alpha_1,\dots,\alpha_r\}$ be a finite collection of
finite subsets of ${\bf S}(A)$ such that the following two conditions
are satisfied:
$$\cap_{k=1}^r\alpha_k\ne\emptyset.$$

For some $m\in\{1,\dots,r-1\}$ and some finite subset
$\alpha\subset{\bf S}(A)$ we have $\alpha\supset\alpha_k$ for
$k\in\{1,\dots,m\}$ and $\alpha\subset\alpha_k$ for
$k\in\{m+1,\dots,r\}.$

Let $n\in{\bf N}$ and $\s\in\cap_{k=1}^r\alpha_k$. We define
$E(\{\alpha_1,\dots,\alpha_r\},\s,n,m)$ as the set of finitely
non-zero vectors of $S(l_1(\Psi))$ satisfying the following
three conditions.

(a) Their supports are contained in the union of the following
two sets:
$${\rm supp}\s\times\{\alpha_1,\dots,\alpha_m\}\times\{1,\dots,n\},$$
$${\rm supp}\s\times\{\alpha_{m+1},\dots,\alpha_r\}\times\{n+1,\dots\}.$$

(b) The signs of all non-zero coordinates corresponding to the triples
$(\gamma,\alpha_k,j)$ of the first set are equal to $\s(\gamma)$.

(c) The signs of all non-zero coordinates corresponding to the triples
$(\gamma,\alpha_k,j)$ of the second set are equal to $-\s(\gamma)$.

The union of all sets $E(\{\alpha_1,\dots,\alpha_r\},\s,n,m)$ constructed
in the described way we shall denote by $E$. It is clear that $E$
satisfies Conditions 1 and 2 of Definition 1.

Let $X$ be a Banach space. Let us show that the following statements
are equivalent.

(a) $X\notin P^{\rm co}$.

(b) $h(X,E)>0.$

(c) $h(X,E)\ge \delta/3.$

It is clear that we need to prove the implications
(a)$\Rightarrow$(c) and (b)$\Rightarrow$(a) only.

Let us start with the implication (a)$\Rightarrow$(c). Let
$0<\theta<1$. By Lemma 1 we can find $u\in S(l_\infty(\Gamma,X^{**}/X))$
such that for every $\s\in{\bf S}(A)$ there exists
$f_\s\in(X^{**}/X)^*=X^\bot\subset X^{***}$ such that
$||f_\s||<1$ and for every $\gamma\in{\rm supp}\s$ we have
$|f_\s(u(\gamma))|>\delta\theta$ and
sign$f_\s(u(\gamma))=\s(\gamma).$ It is clear that all
conditions above are satisfied for some
$u\in l_\infty(\Gamma, X^{**}/X)$ satisfying the inequality
$||u||<1$. So we shall suppose that $||u||<1$.

Let $x^{**}(\gamma)\in X^{**}\ (\gamma\in\Gamma)$ be such that
the image of $x^{**}(\gamma)$ under the quotient mapping
$X^{**}\to X^{**}/X$ coincides with $u(\gamma)$ and
$||x^{**}(\gamma)||<1$.

Now we define $x\in S(l_\infty(\Psi, X))$ for which
$\inf_{e\in E}||e(x)||$ is ``large''. This function will be defined
in stages. At first we define this function for triples
$(\gamma,\alpha,k)$ with one-element set $\alpha$, then for triples
with two-element set $\alpha$, etc.

It is clear that it is sufficient to define $x(\gamma,\alpha,k)$
only for those triples for which
$\gamma\in\cup_{\s\in\alpha}{\rm supp}\s$. (Since for other
triples $(\gamma,\alpha,k)$ the numbers $e(\gamma,\alpha,k)$ are
equal to zero for every $e\in E$ and so $x(\gamma,\alpha,k)$ can
be defined arbitrarily.)

In the course of construction we shall also define functionals
$f(\s,\alpha,k)$, where $\alpha$ is a finite subset of
${\bf S}(A)$, $\s\in\alpha$ and $k\in{\bf N}$.

For one-element sets $\alpha$ we define $x(\gamma,\alpha,k)$ by induction
on $k$.

1. Let $\alpha=\{\s\}.$

1.1. Define $f(\s,\{\s\},1)$ as a reflection of $f_\s$ with respect
to $\{x^{**}(\gamma):\gamma\in{\rm supp}\s\}.$ Let $x(\gamma,\{\s\},1)$
be reflections of $\{x^{**}(\gamma):\gamma\in{\rm supp}\s\}$ with
respect to $\{f(\s,\{\s\},1)\}$.

\dots

1.k. Define $f(\s,\{\s\},k)$ as a reflection of $f_\s$ with respect
to
$$\{x^{**}(\gamma):\gamma\in{\rm supp}\s\}\cup
\{x(\gamma,\{\s\},j):\gamma\in{\rm supp}\s, j=1,\dots,k-1\}.$$

Let $x(\gamma,\{\s\},k)$ be reflections of
$\{x^{**}(\gamma):\gamma\in{\rm supp}\s\}$ with respect to
$$\{f(\s,\{\s\},j):\ j=1,\dots,k\}.$$

\dots

Let us suppose that we have already constructed $x(\gamma,\alpha,k)$
for all $(m-1)$-element subsets $\alpha\subset{\bf S}(A)$.

m. Let $\alpha=\{\s_1,\dots,\s_m\}$. By supp$\alpha$ we shall denote
$\cup_{i=1}^m{\rm supp}\s_m$.

m.1. Define $\{f(\s_i,\alpha,1)\}$ as a reflections of $f_{\s_i}$ with
respect to $\{x^{**}(\gamma):\gamma\in{\rm supp}\alpha\}$. Let
$\{x(\gamma,\alpha,1):\ \gamma\in{\rm supp}\alpha\}$ be
reflections of $\{x^{**}(\gamma):\ \gamma\in{\rm supp}\alpha\}$
with respect to
$\{f(\s,\beta,1):\ \beta\subset\alpha,\ \s\in\beta\}.$

\dots

m.k. Define $\{f(\s_i,\alpha,k)\}$ as a reflections of
$f_{\s_i}$ with respect to
$$\{x^{**}(\gamma):\ \gamma\in{\rm supp}\alpha\}\cup
\{x(\gamma,\beta,j):\ \beta\subset\alpha,\
\gamma\in{\rm supp}\beta,\ j=1,\dots,k-1\}.$$

Let $\{x(\gamma,\alpha,k)\}$ be reflections of
$\{x^{**}(\gamma):\ \gamma\in{\rm supp}\alpha\}$ with respect
to $\{f(\s,\beta,j):\ \beta\subset\alpha,\ \s\in\beta,\
j=1,\dots,k\}$.

\dots

Let us show that
$$(\forall e\in E)(||e(x)||\ge\delta\theta/3).$$

Let $e\in E(\{\alpha_1,\dots,\alpha_r\},\s,n,m),$
where $\s\in\cap_{i=1}^r\alpha_i$. We have
$$e(x)=\sum_{\gamma\in{\rm supp}\s}(\sum_{i=1}^m\sum_{k=1}^n
e(\gamma,\alpha_i,k)x(\gamma,\alpha_i,k)+
\sum_{i=m+1}^r\sum_{k=n+1}^\infty e(\gamma,\alpha_i,k)
x(\gamma,\alpha_i,k)).$$
Therefore we have
$$f(\s,\{\s\},1)(e(x))=
\sum_{\gamma\in{\rm supp}\s}(\sum_{i=1}^m\sum_{k=1}^ne(\gamma,\alpha_i,k)
x^{**}(\gamma)f(\s,\{\s\},1)+$$
$$\sum_{i=m+1}^r\sum_{k=n+1}^\infty e(\gamma,\alpha_i,k)x^{**}(\gamma)
f(\s,\{\s\},1))=$$
$$\sum_{\gamma\in{\rm supp}\s}(\sum_{i=1}^m\sum_{k=1}^n|e(\gamma,\alpha_i,k)|
|f_\s(x^{**}(\gamma))|-
\sum_{i=m+1}^r\sum_{k=n+1}^\infty |e(\gamma,\alpha_i,k)|
|f_\s(x^{**}(\gamma))|)$$

According to definition of $m$ there exists
$\alpha\subset{\bf S}(A)$ such that
$\cup_{i=1}^m\alpha_i\subset\alpha$ and $\alpha\subset\alpha_i$
when $i\ge m+1$. We have
$$f(\s,\alpha,n+1)(e(x))=
\sum_{\gamma\in{\rm supp}\s}(\sum_{i=1}^m\sum_{k=1}^ne(\gamma,\alpha_i,k)
f(\s,\alpha,n+1)x(\gamma,\alpha_i,k)+$$
$$\sum_{i=m+1}^r\sum_{k=n+1}^\infty e(\gamma,\alpha_i,k)
f(\s,\alpha,n+1)x(\gamma,\alpha_i,k)).$$

Since $\alpha_i\subset\alpha$ if $i\le m$ and
$f_\s\in X^\bot$ then the first sum equals zero. Since
$\alpha\subset\alpha_i$ when $i>m$, the second sum can be
rewritten in the following way:
$$\sum_{m+1}^r\sum_{k=n+1}^\infty e(\gamma,\alpha_i,k)
x^{**}(\gamma)f(\s,\alpha,n+1)=
-\sum_{i=m+1}^n\sum_{k=n+1}^\infty|e(\gamma,\alpha_i,k)|
|f_\s(x^{**}(\gamma))|.$$
Hence
$$2|f(\s,\alpha,n+1)(e(x))|+f(\s,\{\s\},1)(e(x))=$$
$$\sum_\gamma\sum_i\sum_k|e(\gamma,\alpha_i,k)||f_\s(x^{**}(\gamma))|>
\delta\theta.$$
Therefore $||e(x)||\ge\delta\theta/3$. Since $\theta$ is arbitrary,
the implication (a)$\Rightarrow$(c) is proved.

Now we turn to the implication (b)$\Rightarrow$(a). Let
$h(X,E)=\tau>0$. Let $\theta$ be arbitrary number from the open
interval (0,1) and let $x\in S(l_\infty(\Psi,X))$ be such that
$\inf_{e\in E}||e(x)||>\tau\theta$.

The set of pairs $(\alpha, k)$, where $\alpha$ is a finite subset of
${\bf S}(A)$ and $k\in{\bf N}$ is a directed set with respect to the
natural order. Let $U$ be some ultrafilter which majorize the filter
induced by this order. Let
$$x^{**}(\gamma)=w^*-\lim_Ux(\gamma,\alpha,k)\ (\gamma\in\Gamma).$$
Thus we defined $x^{**}\in B(l_\infty(\Gamma,X^{**}))$. It is clear
that in order to finish the proof of Theorem 2 it is sufficient to
prove that
$$(\forall a\in A)(||\varphi\sum_\gamma a(\gamma)x^{**}(\gamma)||\ge
\tau\theta/2,$$
where $\varphi$ is the quotient mapping $\varphi:X^{**}\to X^{**}/X$.

Assume the contrary, i.e. there exists $a\in A$ and $x\in X$ such that
$$||\sum_\gamma a(\gamma)x^{**}(\gamma)-x||<\tau\theta/2.$$
By the definition of $x^{**}(\gamma)$ for every finite subset
$\alpha_0\subset{\bf S}(A)$ and every $k_0\in{\bf N}$ we have
$$\sum_\gamma a(\gamma)x^{**}(\gamma)\in w^*-{\rm cl}
\{\sum_\gamma a(\gamma)x(\gamma,\alpha,k):\
\alpha\supset\alpha_0,\ k\ge k_0\}.$$

By Proposition 6.30 of [O1] we can find a collection
$\{\alpha_1,\dots,\alpha_m\}$ of finite subsets of
${\bf S}(A)$, real numbers $\tau_1,\dots,\tau_m\ge 0$ with
$\sum_{i=1}^m\tau_i=1$ and numbers $k_1,\dots,k_m\in{\bf N}$ such that
sign$a\in\cap_{i=1}^m\alpha_m$ and
$$||\sum_\gamma a(\gamma)x^{**}(\gamma)-
\sum_{i=1}^m\tau_i\sum_\gamma a(\gamma)x(\gamma,\alpha_i,k_i)||<
\tau\theta.$$

Using Proposition 6.30 of [O1] once more we find a collection
$\{\alpha_{m+1},\dots,\alpha_r\}$ of finite subsets of
${\bf S}(A)$, satisfying the condition
$\cap_{i=m+1}^r\alpha_i\supset\cup_{i=1}^m\alpha_i$, real numbers
$\tau_{m+1},\dots,\tau_r\ge0$ with $\sum_{i=m+1}^r\tau_i=1$
and numbers $k_{m+1},\dots,k_r\in{\bf N}$ such that
$$\max_{1\le i\le m}k_i<\min_{m+1\le i\le r}k_i$$
and
$$||\sum_\gamma a(\gamma)x^{**}(\gamma)-
\sum_{i=m+1}^r\tau_i\sum_\gamma a(\gamma)x(\gamma,\alpha_i,k_i)||<
\tau\theta.$$
Hence
$$||\sum_{i=1}^m\tau_i\sum_\gamma a(\gamma)x(\gamma,\alpha_i,k_i)-
\sum_{i=m+1}^r\tau_i\sum_\gamma a(\gamma)x(\gamma,\alpha_i,k_i)||<
2\tau\theta.$$
It remains to note that the vector in the left-hand side of this
inequality has the form $2e(x)$ for some $e\in E$. We arrived at a
contradiction. The theorem is proved.$\Box$

\centerline{\bf References}

[AAG] J.Alvarez,  T.Alvarez  and  M.Gonzalez, The  gap   between
subspaces and perturbation of non semi-\-Fred\-holm operators,
Bull. Austral. Math. Soc.
45
(1992)
no. 3, pp. 369--376

[B1]
B.Beauzamy, Introduction  to  Banach  Spaces   and   their
Geometry, Amsterdam,
North-Holland Publishing Company,
1982

[B2]
B.Beauzamy, Banach-Saks  properties  and  spreading  models,
Math.\- Scand.
44
(1979)
no. 2,
pp. 357--384

[J] R.Janz,
Pertubation of Banach  spaces, preprint,  Konstanz,
1987

[K] M.I.Kadets,
Note on the gap between subspaces
(Russian),
Funkts. Anal. Prilozhen.
9
(1975)
no. 2,
pp. 73--74

[KKM]
M.G.Krein, M.A.Krasnoselskii and D.P.Milman,
On  the  defect
numbers of linear operators  in  a  Banach  space  and  on  some
geometric questions (Russian),
Sbornik Trudov Inst.  Matem.  AN  Ukrainian
SSR
11
(1948)
pp. 97--112

[O1] M.I.Ostrovskii,
Topologies on the set of  all  subspaces  of  a
Banach space and  related  questions  of  Banach  space  geometry,
Quaestiones Math. (to appear)

[O2] M.I.Ostrovskii,
The relation of  the  gaps  of  subspaces  to
dimension and other isomorphic invariants (Russian), Akad.  Nauk  Armyan.
SSR Dokl. 79
(1984)
no. 2
pp. 51--53

[P]
A.Pe\l czy\'nski
All separable Banach spaces admit for every
$\varepsilon >0$
fundamental total  and  bounded  by $1+\varepsilon $  biorthogonal
sequences,
Studia Math.
55
(1976)
pp. 295--304

\end{document}